\documentclass[11pt]{article}

\usepackage[T2A]{fontenc}
\usepackage[cp1251]{inputenc}
\usepackage{amsfonts}
\usepackage{eufrak}
\usepackage{amssymb}
\usepackage{amsmath}



\begin{document}

\newtheorem{theorem}{Theorem}
\newtheorem{lemma}{Lemma}
\newtheorem{corollary}{Corollary}
\newtheorem{definition}{Опеределение}
\newtheorem{proposition}{Proposition}
\newtheorem{remark}{Remark}

\noindent
\verb"The final version of the article will be published in "

\noindent
\verb"Journal of Difference Equations and Applications."

\bigskip

\centerline{\textbf{On  analogues of   C.R.~Rao's theorems}}

\centerline{\textbf{for locally compact Abelian groups}}

\bigskip

\centerline{\bf G.M. Feldman}

\bigskip

 \makebox[20mm]{ }\parbox{125mm}{ \small Let $\xi_1$, $\xi_2$, $\xi_3$ be independent random variables with nonvanishing characteristic functions, and $a_j$, $b_j$ be real numbers such that $a_i/b_i\ne a_j/b_j$ for $i\ne j$. Let $L_1=a_1\xi_1+a_2\xi_2+a_3\xi_3$,  $L_2=b_1\xi_1+b_2\xi_2+b_3\xi_3$. By C.R.~Rao's theorem   the distribution of the random vector  $(L_1, L_2)$ determines the distributions  of the random variables $\xi_j$ up to a  change of location.  We prove an analogue of this theorem for   independent random variables with values in a locally compact Abelian group. We also prove an analogue of similar C.R.~Rao's theorem for independent random variables with values in an \text{\boldmath $a$}-adic solenoid. In so doing coefficients of linear forms are continuous endomorphisms  of the group.}

\bigskip


\bigskip


\bigskip

\centerline{\bf 1. Introduction}

\bigskip

In the article \cite{CRR} C.R.~Rao proved the following theorem.

\medskip

{\bf Theorem A}. {\it  Let $\xi_1$, $\xi_2$, $\xi_3$ be independent random variables with nonvanishing characteristic functions, and $a_j$, $b_j$ be real numbers such that $a_i/b_i\ne a_j/b_j$ for $i\ne j$. Let $L_1=a_1\xi_1+a_2\xi_2+a_3\xi_3$,  $L_2=b_1\xi_1+b_2\xi_2+b_3\xi_3$. Then the distribution of the random vector  $(L_1, L_2)$ determines the distributions  of the random variables $\xi_j$ up to a  change of location.}

\medskip

The C.R.~Rao theorem is a generalisation of the following result by I.~Kotlarski \cite{K}.

\medskip

{\bf Theorem B}. {\it      Let $\xi_1$, $\xi_2$, $\xi_3$ be independent random variables with nonvanishing characteristic functions. Let $L_1=\xi_1-\xi_3$,  $L_2=\xi_2-\xi_3$. Then the distribution of the random vector  $(L_1, L_2)$ determines the distributions  of the random variables $\xi_j$ up to a  change of location.}

\medskip

The aim of this note is to prove   an analogue of the C.R.~Rao theorem for  linear forms of independent random variables with values in a locally compact Abelian group. We also prove an analogue for independent random variables with values in an \text{\boldmath $a$}-adic solenoid of similar C.R.~Rao's theorem.
In so doing coefficients of linear forms are continuous endomorphisms  of the group.

Let $X$ be a second countable locally compact Abelian group. We will consider only such groups, without mentioning it specifically.   Denote by $Y$ the character
group of the group $X$, and by  $(x,y)$ the value of a character $y \in Y$ at an element $x \in X$.
Let
 $\alpha:X \mapsto X $
be a continuous endomorphism.
The adjoint endomorphism $\tilde\alpha: Y \mapsto Y $
is defined by the formula $(\alpha x,
y)=(x , \tilde\alpha y)$ for all $x\in X$, $y\in
Y$.  If $G$ is a closed subgroup of $X$, denote by
 $A(Y, G) = \{y \in Y: (x, y) = 1$ \mbox{ for all } $x \in G \}$
its annihilator. Denote by   $\mathbb{T}$ the circle group (the one dimensional torus) and by $\mathbb{Z}$ the group of integers.

Let $f(y)$ be a function on the group    $Y$,   and let $h$ be an arbitrary element of $Y$. Denote by   $\Delta_h$   the finite difference operator
$$
\Delta_h f(y)=f(y+h)-f(y).
$$
A  function $f(y)$ on the group $Y$ is called a
polynomial if for some  $n$ it satisfies the equation
$$\Delta_{h}^{n+1}f(y)=0, \ y,h \in Y.$$
The minimal
$n$ for which this equality takes place is called the degree of $f(y)$.

Denote by ${\rm M}^1(X)$ the
convolution semigroup of probability  distributions  on the group $X$. Let
${\mu\in {\rm M}^1(X)}$.  Denote by
$$
\hat\mu(y) =
\int_{X}(x, y)d \mu(x), \ y\in Y,$$ the characteristic function
of  the distribution $\mu$. Denote by   $E_x$  the degenerate distribution
 concentrated at a point $x\in X$.

A distribution  $\gamma\in {\rm M}^1(X)$ is called Gaussian
(\cite[Chapter IV]{Pa}),  
if its characteristic function is represented in the form
\begin{equation}\label{f1}
\hat\gamma(y)= (x,y)\exp\{-\varphi(y)\}, \  y\in Y,
\end{equation}
where $x \in X$, and $\varphi(y)$ is a continuous nonnegative function
on the group $Y$
 satisfying the equation
 \begin{equation}\label{f2}
\varphi(u + v) + \varphi(u
- v) = 2[\varphi(u) + \varphi(v)], \ u,  v \in
Y.
\end{equation}

\bigskip

\centerline{\bf 2. Proof of the main theorem}

\bigskip

The following result can be considered as an analogue of Theorem A for locally compact Abelian groups.

\begin{theorem}\label{th1}.
\textit{Let $X$ be a locally compact Abelian group, $\xi_1$, $\xi_2$, $\xi_3$ be independent random variables with values in the group   $X$ with nonvanishing characteristic functions, and $b_j$ be continuous endomorphisms of $X$. Let $L_2=b_1\xi_1+b_2\xi_2+b_3\xi_3$. Assume that either}

{$(I)$ $L_1=\xi_1+\xi_2+\xi_3$ and $b_j$ satisfy the conditions
\begin{equation}\label{08_03_1}
{\rm Ker} (b_1-b_2)={\rm Ker} (b_1-b_3)={\rm Ker} (b_2-b_3)=\{0\},
\end{equation}
or}

{$(II)$ $L_1=\xi_1+\xi_2$ and $b_j$ satisfy the conditions
\begin{equation}\label{07_03_1a}
{\rm Ker} (b_1-b_2)={\rm Ker}~b_3=\{0\}.
\end{equation}
 Then the distribution of the random vector  $(L_1, L_2)$ determines the distributions  of the random variables $\xi_j$ up to a  shift.}
\end{theorem}

For the convenience of references  we formulate the following assertion in the form of a lemma.
\begin{lemma}\label{l3} {\rm (\cite[(24.41)]{Hewitt-Ross})}. {\it  Let $X$ be a locally compact Abelian group,  $Y$ be its character group. Let $a$ be a continuous endomorphism of the group $X$.  The set ${\tilde a (Y)}$ is dense in $Y$
if an only if
${\rm Ker}~a=\{0\}$.
}
\end{lemma}

{\bf Proof of Theorem \ref{th1}}.  Let $\eta_1$, $\eta_2$, $\eta_3$ be independent random variables with values in the group   $X$ with nonvanishing characteristic functions. Denote by $\mu_j$ the distribution of the random variable   $\xi_j$, and by $\nu_j$
the distribution of the random variable $\eta_j$.  Put $M_2=b_1\eta_1+b_2\eta_2+b_3\eta_3$,
 $f_j(y)=\hat\nu_j(y)/\hat\mu_j(y)$, $j=1, 2, 3$.

Assume that $(I)$ holds. Put $M_1=\eta_1+\eta_2+\eta_3$. The characteristic function of the random vector   $(L_1, L_2)$  is of the form
\begin{equation}\label{05_03_1}
{\bf E}[(L_1, u)(L_2, v)]={\bf E}[(\xi_1+\xi_2+\xi_3, u)(b_1\xi_1+b_2\xi_2+b_3\xi_3, v)]$$$$
={\bf E}[(\xi_1,u+\tilde b_1v)(\xi_2,u+\tilde b_2v)(\xi_3,u+\tilde b_3v)]
=\hat\mu_1(u+\tilde b_1v)\hat\mu_2(u+\tilde b_2v)\hat\mu_3(u+\tilde b_3v), \ u, v\in Y.
\end{equation}
Assume that the random vectors  $(L_1, L_2)$ and $(M_1, M_2)$ are identically distributed, and hence, they have the same  characteristic functions. Then it follows from (\ref{05_03_1}) that
\begin{equation}\label{05_03_2}
\hat\mu_1(u+\tilde b_1v)\hat\mu_2(u+\tilde b_2v)\hat\mu_3(u+\tilde b_3v)=\hat\nu_1(u+\tilde b_1v)\hat\nu_2(u+\tilde b_2v)\hat\nu_3(u+\tilde b_3v), \ u, v\in Y.
\end{equation}
 We find from    (\ref{05_03_2}) that the functions $f_j(y)$ satisfy the equation
\begin{equation}\label{05_03_3}
f_1(u+\tilde b_1v)f_2(u+\tilde b_2v)f_3(u+\tilde b_3v)=1, \ u, v\in Y.
\end{equation}

We will prove that there exist elements   $x_j\in X$  such that $f_j(y)=(x_j, y)$, $y\in Y$, $j=1, 2, 3$. To prove we use the finite difference method.
Let $k_1$ be an arbitrary element of the group
$Y$. Put $h_1=-\tilde{b}_1k_1$. Then
$h_1+\tilde{b}_1k_1=0$. Substitute $u+h_1$ for $u$ and
$v+k_1$ for $v$  in equation (\ref{05_03_3}).
Dividing the resulting equation by equation (\ref{05_03_3}), we obtain
\begin{equation}\label{05_03_4}
{f_2(u+\tilde b_2v+(\tilde b_2-\tilde b_1)k_1)f_3(u+\tilde b_3v+(\tilde b_3-\tilde b_1)k_1)\over
f_2(u+\tilde b_2v)f_3(u+\tilde b_3v)}=1, \ u, v\in Y.
\end{equation}
Let $k_2$ be an arbitrary element of the group
$Y$. Put $h_2=-\tilde{b}_2k_2$. Then
$h_2+\tilde{b}_2k_2=0$. Substitute $u+h_2$ for $u$ and
$v+k_2$ for $v$  in equation (\ref{05_03_4}).
Dividing the resulting equation by equation (\ref{05_03_4}), we find
\begin{equation}\label{05_03_5}
{f_3(u+\tilde b_3v+(\tilde b_3-\tilde b_1)k_1+(\tilde b_3-\tilde b_2)k_2)f_3(u+\tilde b_3v)\over
f_3(u+\tilde b_3v+(\tilde b_3-\tilde b_1)k_1)f_3(u+\tilde b_3v+(\tilde b_3-\tilde b_2)k_2)}=1, \ u, v\in Y.
\end{equation}
 Substitute $u=v=0$, $(\tilde b_3-\tilde b_1)k_1=k$, $(\tilde b_3-\tilde b_2)k_2=l$ in equation (\ref{05_03_5}). Note that since
  \begin{equation}\label{27_07_1}
{\rm Ker} (b_1-b_3)={\rm Ker} (b_2-b_3)=\{0\},
\end{equation}
  by Lemma \ref{l3} the sets ${(\tilde b_3-\tilde b_1)(Y)}$ and ${(\tilde b_3-\tilde b_2)(Y)}$ are dense in $Y$. Since $k_1$ and  $k_2$ are arbitrary elements of the group   $Y$, this implies that the function $f_3(y)$ satisfies the equation
$$
{f_3(k +l)\over
f_3(k)f_3(l)}=1, \ k, l\in Y,
$$
and hence,
\begin{equation}\label{10_03_1}
{f_3(k+l)=f_3(k)f_3(l)}, \ k,  l\in Y.
\end{equation}
Taking into account that $f_3(-y)=\overline{f_3(y)}$, $y\in Y$, we find from equation (\ref{10_03_1}) that $|f_3(y)|=1$, $y\in Y$. Thus, we proved that the function $f_3(y)$ is a character of the group $Y$. By the Pontryagin duality theorem there is an element $x_3\in X$ such that $f_3(y)=(x_3, y).$ Note that this proof  uses  only conditions (\ref{27_07_1}), but not  ${\rm Ker} (b_1-b_2)=\{0\}$ (compare below with Remark \ref{re3}, where   we discuss the necessity of conditions (\ref{08_03_1}) and  (\ref{07_03_1a}) for     Theorem  \ref{th1} to be valid).

Arguing in a similar way we get that there exist elements $x_j\in X$ such that $f_j(y)=(x_j, y)$, $y\in Y$, $j=1, 2$.
Thus, $\hat\nu_j(y)=\hat\mu_j(y)(x_j, y)$, $y\in Y$, and hence, $\nu_j=\mu_j*E_{x_j}$, $j=1, 2, 3$. Thus, in the case when $(I)$ holds, the theorem is proved.

Assume that $(II)$ holds. Put $M_1=\eta_1+\eta_2$.
Arguing as in the proof of the theorem in the case when $(I)$ holds, we get that the characteristic function of the random vector $(L_1, L_2)$  is of the form
\begin{equation}\label{08_03_6}
{\bf E}[(L_1, u)(L_2, v)]=\hat\mu_1(u+\tilde b_1v)\hat\mu_2(u+\tilde b_2v)\hat\mu_3(\tilde b_3v), \ u, v\in Y,
\end{equation}
and the identical distributiveness of the random vectors    $(L_1, L_2)$ and $(M_1, M_2)$ implies that the functions $f_j(y)$ satisfy the equation
\begin{equation}\label{06_03_1}
f_1(u+\tilde b_1v)f_2(u+\tilde b_2v)f_3(\tilde b_3v)=1, \ u, v\in Y.
\end{equation}
The theorem will be proved if we prove that the functions $f_j(y)$, $j=1, 2, 3$,  are characters of the group $Y$.  First verify that the function  $f_3(y)$ is a character of the group $Y$. We follow the scheme of the proof of the theorem in the case when $(I)$ holds.   Eliminating from equation  (\ref{06_03_1}) the function $f_1(y)$, we obtain
\begin{equation}\label{08_03_4}
{f_2(u+\tilde b_2v+(\tilde b_2-\tilde b_1)k_1)f_3(\tilde b_3v+\tilde b_3k_1)\over
f_2(u+\tilde b_2v)f_3(\tilde b_3v)}=1, \ u, v\in Y.
\end{equation}
Eliminating from equation  (\ref{08_03_4}) the function$f_2(y)$, we get
\begin{equation}\label{08_03_5}
{f_3(\tilde b_3v+\tilde b_3k_1+\tilde b_3k_2)f_3(\tilde b_3v)\over
f_3(\tilde b_3v+\tilde b_3k_1)f_3(\tilde b_3v+\tilde b_3k_2)}=1, \ u, v\in Y.
\end{equation}
Putting $v=0$ in equation  (\ref{08_03_5}), using Lemma \ref{l3} and the fact that ${\rm Ker}~b_3=\{0\}$, we make sure that the function $f_3(y)$ is a character of the group $Y$.

Verify now that the functions $f_j(y)$, $j=1, 2$, are characters of the group  $Y$. Let $h$ be an arbitrary element of the group
$Y$.   Substitute $u+h$ for $u$  in equation (\ref{08_03_4}).
Dividing the resulting equation into equation (\ref{08_03_4}), we obtain
$$
{f_2(u+\tilde b_2v+(\tilde b_2-\tilde b_1)k_1+h)f_2(u+\tilde b_2v)\over
f_2(u+\tilde b_2v+h)f_2(u+\tilde b_2v+(\tilde b_2-\tilde b_1)k_1)}=1, \ u, v\in Y.
$$
Putting here $u=v=0$, using Lemma \ref{l3} and the fact that ${\rm Ker} (b_2-b_1)=\{0\}$,
we make sure that the function $f_2(y)$ is a character of the group $Y$.
 Arguing in a similar way we prove that the function $f_1(y)$ is a character of the group $Y$. The theorem is completely proved.
$\Box$
\begin{remark}\label{re1}.
{\rm Suppose that  in Theorem \ref{th1} $L_1=\xi_1+\xi_2$ and $L_2=\xi_2+\xi_3$. Then, obviously, continuous endomorphisms $b_j$ satisfy conditions  (\ref{07_03_1a}), and Theorem \ref{th1} implies, in particular, an   analogue of the I.~Kotlarsi theorem for locally compact Abelian groups (see \cite{PR}).}
\end{remark}
\begin{remark}\label{re2}.
{\rm
It is easy to see that Theorem  \ref{th1}    implies the following theorem proved by    C.R.~Rao in   \cite{CRR}.} \end{remark}

{\bf Theorem C}. {\it  Let   $\xi_1$, $\xi_2$, $\xi_3$ be independent random vectors with values in the space $\mathbb{R}^p$  with nonvanishing characteristic functions. Let $a_j, b_j$ be $(p\times p)$-matrices satisfying the conditions:

$(i)$ $a_j$ is either the zero matrix or a nonsingular matrix, and the only one of the matrices $a_j$  may be the zero matrix;

$(ii)$ $b_j$ is either the zero matrix or a nonsingular matrix, and the only one of the matrices  $b_j$   may be the zero matrix;

$(iii)$ the matrices $a_j$ and $b_j$ can not be both the zero matrices;

$(iv)$ if the matrix $b_i{a_i}^{-1}-b_j{a_j}^{-1}$, $i\ne j$,  is defined, then it is nonsingular.

Let  $L_1=a_1\xi_1+a_2\xi_2+a_3\xi_3$,  $L_2=b_1\xi_1+b_2\xi_2+b_3\xi_3$.  Then the distribution of the pair    $(L_1, L_2)$  determines the distributions  of the random vectors $\xi_j$ up to a  change of location.}

\medskip

We supplement Theorem  \ref{th1}   with the following assertion.
 \begin{proposition}\label{pr1}.
 {\it Let $X$ be a locally compact Abelian group, $\xi_1$ and  $\xi_2$ be independent random variables with values in the group   $X$ with nonvanishing characteristic functions. Let $b_1$ and $b_2$ be continuous endomorphisms of the group $X$  such that
 \begin{equation}\label{08_03_2}
{\rm Ker} (b_1-b_2)=\{0\}.
\end{equation}
Let $L_1=\xi_1+\xi_2$,  $L_2=b_1\xi_1+b_2\xi_2$. Then the distribution of the random vector  $(L_1, L_2)$ uniquely determines the distributions  of the random variables  $\xi_1$ and $\xi_2$}.
\end{proposition}
  {\bf Proof}. Retaining the notation used in the proof of  Theorem \ref{th1},   we get that the  functions $f_j(y)$ satisfy the equation
$$
f_1(u+\tilde b_1v)f_2(u+\tilde b_2v)=1, \ u, v\in Y.
$$
Putting here $u=-\tilde b_2y$, $v=y$, we obtain
\begin{equation}\label{07_03_1}
f_1((\tilde b_1-\tilde b_2)y)=1, \ y \in Y.
\end{equation}
 Since ${\rm Ker} (b_1-b_2)=\{0\}$, by Lemma \ref{l3} the set ${(\tilde b_1-\tilde b_2)(Y)}$   is dense in $Y$. Then equation  (\ref{07_03_1})
implies that
$f_1(y)=1$, $y\in Y$. Thus, $\hat\nu_1(y)=\hat\mu_1(y)$, $y\in Y$, i.e. $\nu_1=\mu_1$. Arguing in a similar way we get that $\nu_2=\mu_2$. $\Box$
\begin{remark}\label{re3}.
{\rm Theorem  \ref{th1} and   Proposition \ref{pr1} fail  if we omit conditions   (\ref{08_03_1}), (\ref{07_03_1a}) and (\ref{08_03_2}) respectively. Construct the corresponding examples.  Recall the definition of a  Poisson distribution on the group $X$. Let  $x_0\in X$, $\lambda>0$.  Put
$$
e(\lambda E_{x_0})=e^{-\lambda}\left(E_0+\lambda E_{x_0}+\lambda^2 E_{2x_0}/2!+
\dots+\lambda^nE_{nx_0}/n!+\dots\right).
$$
 The characteristic function of the   Poisson distribution  $e(\lambda E_{x_0})$ is of the form
\begin{equation}\label{08_03_3}
\widehat{e(\lambda E_{x_0})}(y)=\exp\left\{\lambda((x_0, y)-1)\right\}, \ y\in Y.
\end{equation}
Let $b_j$, $j=1, 2, 3$, be continuous endomorphisms of the group $X$. Assume that conditions  (\ref{08_03_1}) are not satisfied. Assume for definiteness that  ${\rm Ker} (b_1-b_2)\ne\{0\}$. Take $x_0\in {\rm Ker} (b_1-b_2)$, $x_0\ne 0$, and a number $a>0$. Consider the distributions $\mu_1=\mu_2=e(2aE_{x_0})$, $\nu_1=e(aE_{x_0})$, $\nu_2=e(3aE_{x_0})$, $\mu_3=\nu_3=\mu$, where $\mu$ is an arbitrary distribution on $X$. Let    $\xi_1$, $\xi_2$, $\xi_3$ be independent random variables with values in the group   $X$ and distributions $\mu_j$, and $\eta_1$, $\eta_2$, $\eta_3$ be independent random variables with values in the group   $X$ and distributions $\nu_j$. It is obvious that $\nu_1$ and $\nu_2$  are not shifts of  $\mu_1$ and $\mu_2$. Put $L_1=\xi_1+\xi_2+\xi_3$,  $L_2=b_1\xi_1+b_2\xi_2+b_3\xi_3$, $M_1=\eta_1+\eta_2+\eta_3$,  $M_2=b_1\eta_1+b_2\eta_2+b_3\eta_3$. Since $x_0\in {\rm Ker} (b_1-b_2)$, we have $b_1x_0=b_2x_0=\tilde x$.  It follows from  (\ref{05_03_1}) and (\ref{08_03_3}) that
$${\bf E}[(L_1, u)(L_2, v)]=\exp\{2a((x_0, u+\tilde b_1v)-1)\}\exp\{2a((x_0, u+\tilde b_2v)-1)\}\hat\mu(u+\tilde b_3v)$$$$=e^{-4a}\exp\{4a(x_0, u)(\tilde x, v)\}\hat\mu(u+\tilde b_3v), \ u, v\in Y,$$
 $${\bf E}[(M_1, u)(M_2, v)]=\exp\{a((x_0, u+\tilde b_1v)-1)\}\exp\{3a((x_0, u+\tilde b_2v)-1)\}\hat\mu(u+\tilde b_3v)$$$$=e^{-4a}\exp\{4a(x_0, u)(\tilde x, v)\}\hat\mu(u+\tilde b_3v), \ u, v\in Y.$$
 We see that the characteristic functions  of the random vectors $(L_1, L_2)$ and $(M_1, M_2)$  coincide. Hence, the distributions of the random vectors $(L_1, L_2)$ and $(M_1, M_2)$ also coincide. Thus, if $L_1=\xi_1+\xi_2+\xi_3$ and conditions  (\ref{08_03_1}) are not satisfied, then Theorem \ref{th1} fails. Note also that if ${\rm Ker} (b_1-b_2)\ne\{0\}$  and   ${\rm Ker} (b_1-b_3)\ne\{0\}$, then using the above argument, it is not difficult to construct   independent random variables $\xi_1$, $\xi_2$, $\xi_3$ with values in the group   $X$ and distributions $\mu_j$ and independent random variables  $\eta_1$, $\eta_2$, $\eta_3$ with values in the group   $X$ and distributions  $\nu_j$, such that the distributions of the random vectors $(L_1, L_2)$ and $(M_1, M_2)$  coincide, but none of the distributions  $\nu_j$ is a shift of the distribution $\mu_j$.

 Arguing in a similar way we convince ourselves that Theorem \ref{th1}  fails in the case when $L_1=\xi_1+\xi_2$  if we omit the condition ${\rm Ker} (b_1-b_2)=\{0\}$. For Proposition \ref{pr1} the reasoning is the same.

Assume now that $G={\rm Ker}~b_3\ne\{0\}$. Let $\mu_3$ and $\nu_3$ be arbitrary distributions on the group $X$ supported in $G$. It is easy to see that then $\hat\mu_3(y)=\hat\nu_3(y)=1$, $y\in A(Y, G)$. Since $\overline{\tilde a(Y)}=A(Y, {\rm Ker}~a)$ for any continuous endomorphism   $a$ of the group $X$  (\cite[(24.38)]{Hewitt-Ross}), we have $\overline{\tilde b_3(Y)}=A(Y, G)$. Hence
\begin{equation}\label{08_03_7}
\hat\mu_3(\tilde b_3y)=\hat\nu_3(\tilde b_3y)=1, \ y\in Y.
\end{equation}
 Let $\mu_1=\nu_1$ and $\mu_2=\nu_2$ be arbitrary distributions on the group $X$, and assume that $\nu_3$ is not a shift of  $\mu_3$.
Let    $\xi_1$, $\xi_2$, $\xi_3$ be independent random variables with values in the group   $X$ and distributions $\mu_j$, and   $\eta_1$, $\eta_2$, $\eta_3$ be independent random variables with values in the group   $X$ and distributions  $\nu_j$. Put $L_1=\xi_1+\xi_2$,  $L_2=b_1\xi_1+b_2\xi_2+b_3\xi_3$, $M_1=\eta_1+\eta_2$,  $M_2=b_1\eta_1+b_2\eta_2+b_3\eta_3$. It follows from (\ref{08_03_6}) and (\ref{08_03_7}) that the characteristic functions  of the random vectors   $(L_1, L_2)$ and  $(M_1, M_2)$  coincide. Hence, the distributions of the random vectors $(L_1, L_2)$ and $(M_1, M_2)$ also coincide. Thus,  Theorem   \ref{th1}  also fails if  we omit the condition ${\rm Ker}~b_3=\{0\}$.
}\end{remark}

\bigskip

\centerline{\bf 3. C.R.~Rao's theorem for \text{\boldmath $a$}-adic solenoids}

\bigskip

In the article   \cite{CRR}   C.R.~Rao proved the following theorem.

\medskip

{\bf Theorem D}. {\it  Let $\xi_j$, $j=1, 2, 3, 4,$ be independent random variables   with nonvanishing characteristic functions, and  $a_j$, $b_j$ be real numbers such that  $a_i/b_i\ne a_j/b_j$ for $i\ne j$. Let $L_1=a_1\xi_1+a_2\xi_2+a_3\xi_3+a_4\xi_4$,  $L_2=b_1\xi_1+b_2\xi_2+b_3\xi_3+b_4\xi_4$. Then the distribution of the random vector  $(L_1, L_2)$ determines the distributions  of the random variables $\xi_j$ up to  convolution with a Gaussian distribution.}

\medskip

We prove that Theorem D is valid for   \text{\boldmath $a$}-adic solenoids. Recall the definition of an \text{\boldmath $a$}-adic solenoid.
Put \text{\boldmath $a$}= $(a_0, a_1,\dots)$, where all $a_j \in {\mathbb{Z}}$,
$a_j > 1$. Denote by $\Delta_{{\text{\boldmath $a$}}}$
 the group of ${\text{\boldmath $a$}}$-adic integers.  Consider the group
$\mathbb{R}\times\Delta_{{\text{\boldmath $a$}}}$. Denote by
 $B$ the subgroup of the group
 $\mathbb{R}\times\Delta_{{\text{\boldmath $a$}}}$ of the form
$B=\{(n,n\mathbf{u})\}_{n=-\infty}^{\infty}$, where
$\mathbf{u}=(1, 0,\dots,0,\dots)$. The factor-group $\Sigma_{{\text{\boldmath $a$}}}=(\mathbb{R}\times\Delta_{{\text{\boldmath $a$}}})/B$ is called
  an ${\text{\boldmath $a$}}$-{adic solenoid}.   The group $\Sigma_{{\text{\boldmath $a$}}}$ is  compact, connected and has
dimension 1     (\cite[(10.12), (10.13),
(24.28)]{Hewitt-Ross}). The character group of the group  $\Sigma_{{\text{\boldmath $a$}}}$  is topologically isomorphic to a discrete group of the form
$$H_{\text{\boldmath $a$}}= \left\{{m \over a_0a_1 \dots a_n} : \ n = 0, 1,\dots; \ m
\in {\mathbb{Z}} \right\}.
$$
In order not to complicate the notation we will identify $H_{\text{\boldmath $a$}}$ with the character group of the group $\Sigma_{{\text{\boldmath $a$}}}$.

It follows from  (\ref{f1}) and (\ref{f2}) that the characteristic function of a Gaussian distribution   $\gamma$ on  an   \text{\boldmath $a$}-adic solenoid
 $\Sigma_\text{\boldmath $a$}$ is of the form
$$
\hat\gamma(y)=(x, y)\exp\{-\sigma y^2\}, \ y\in H_{\text{\boldmath $a$}},
$$
where $x\in  \Sigma_\text{\boldmath $a$}$, $\sigma \ge 0$.

The following result can be considered as an analogue of Theorem D for \text{\boldmath $a$}-adic solenoids.
\begin{theorem}\label{th3}.
\textit{Let $X=\Sigma_{{\text{\boldmath $a$}}}$ be an \text{\boldmath $a$}-adic solenoid. Let $\xi_j$, $j=1, 2, 3, 4,$ be independent random variables with values in the group   $X$ with nonvanishing characteristic functions, and $b_j$ be continuous endomorphisms of the group $X$.  Let $L_2=b_1\xi_1+b_2\xi_2+b_3\xi_3+b_4\xi_4$. Assume that either}

{$(I)$ $L_1=\xi_1+\xi_2+\xi_3+\xi_4$ and $b_j$ satisfy the conditions
\begin{equation}\label{25_03_1}
{\rm Ker} (b_i-b_j)=\{0\}, \ i\ne j, \  i,j\in\{1, 2, 3, 4\},
\end{equation}
or

$(II)$ $L_1=\xi_1+\xi_2+\xi_3$ and $b_j$ satisfy the conditions
\begin{equation}\label{25_03_2}
{\rm Ker} (b_i-b_j)=\{0\}, \  i\ne j, \ i,j\in\{1, 2, 3\}, \ {\rm Ker}~b_4=\{0\}.
\end{equation}
Then the distribution of the random vector  $(L_1, L_2)$ determines the distributions  of the random variables $\xi_j$ up to  convolution with a Gaussian distribution.}
\end{theorem}

To prove the theorem we need the following lemmas.

\begin{lemma}\label{l1} {\rm (compare with \cite[Lemma 3]{F1})}. {\it  Let $Y$ be a locally compact Abelian group, $\beta_j$ be continuous endomorphisms of the group $Y$ satisfying the conditions:   the sets ${(\beta_i-\beta_j)(Y)}$  are dense in $Y$  for $i\ne j$.  Consider on the group
$Y$ the equation
\begin{equation}\label{09_03_1}
\sum_{j = 1}^{n}  \psi_j(u + \beta_j v ) = B(v),
\ u, v \in Y,
\end{equation}
where $\psi_j(y)$ and  $B(y)$ are continuous functions on $Y$. Then
$\psi_j(y)$ are polynomial on $Y$ of degree at most $n-1$. If $B(y)=0$ for  $y\in Y$,   then $\psi_j(y)$ are of degree at most $n-2$.
 }
\end{lemma}
{\bf Proof}. To prove we use the finite difference method. Let $k_1$ be an arbitrary element of the group
$Y$. Put $h_1=-{\beta_n}k_1$. Then
$h_1+\beta_n k_1=0$. Substitute $u+h_1$ for $u$ and
$v+k_1$ for $v$  in equation (\ref{09_03_1}). Subtracting equation
(\ref{09_03_1})  from the   resulting equation we get
\begin{equation}
\label{09_03_2}
    \sum_{j = 1}^{n-1} \Delta_{l_{1j}}{\psi_j(u + \beta_j v)}
    =\Delta_{k_1}B(v),
\ u,v\in Y,
\end{equation}
where $l_{1j}= h_1+\beta_j k_1=(\beta_j-\beta_n)k_1$, $j=1, 2,\dots, n-1$.  Let
$k_2$
be an arbitrary element of the group
$Y$. Put
$h_2=-{\beta_{n-1}} k_2$. Then $h_2+\beta_{n-1} k_2 =0$. Substitute $u+h_2$ for $u$ and
$v+k_2$ for $v$  in equation   (\ref{09_03_2}).    Subtracting equation (\ref{09_03_2}) from the   resulting equation we obtain
$$
    \sum_{j = 1}^{n-2} \Delta_{l_{2j}} \Delta_{l_{1j}}{\psi_j(u + \beta_j v)}
    =\Delta_{k_2}\Delta_{k_1}B(v),
\ u, v\in Y,
$$
where $l_{2j}= h_2+\beta_j k_2=(\beta_j-\beta_{n-1})k_2$, $j=1, 2,\dots,n-2$.
Arguing in a similar way we arrive at the equation
\begin{equation}\label{09_03_3}
   \Delta_{l_{n-1,1}} \Delta_{l_{n-2,1}}\dots \Delta_{l_{11}}
   {\psi_1(u + \beta_1 v)}
    =\Delta_{k_{n-1}}\Delta_{k_{n-2}}\dots\Delta_{k_1}B(v),
\ u, v\in Y,
\end{equation}
where $k_m$ are arbitrary elements of the group $Y$,
$h_m=-{\beta_{n-m+1}}k_m$,
   $m=1, 2,\dots, n-1$, $l_{mj}=
h_m+\beta_j k_m=(\beta_j-\beta_{n-m+1})k_m$, $j=1, 2,\dots, n-m$.
Let $h$ be an arbitrary element of the group
$Y$. Substitute $u+h$ for $u$  in equation   (\ref{09_03_3}). Subtracting equation (\ref{09_03_3}) from the   resulting equation we obtain
\begin{equation}\label{25_03_3}
   \Delta_{h}\Delta_{l_{n-1,1}} \Delta_{l_{n-2,1}}\dots \Delta_{l_{11}}
   {\psi_1(u + \beta_1 v)}
    =0,
\ u, v\in Y.
\end{equation}
Putting $v=0$ in equation (\ref{25_03_3})  and taking into account that all sets ${(\beta_i-\beta_j)(Y)}$  are dense in $Y$  for $i\ne j$, we obtain that the function $\psi_1(y)$ satisfies the equation
$$\Delta_{h}^{n}\psi_1(y)=0, \ y,h \in Y,$$
i.e. $\psi_1(y)$ is a polynomial on $Y$ of degree at most $n-1$.  If $B(y)=0$ for  $y\in Y$, then in (\ref{09_03_3}) the right-hand side is equal to zero, and the polynomial $\psi_1(y)$ is  of degree at most $n-2$. For the functions  $\psi_j(y)$, $j=2, 3, \dots, n$, we argue in a similar way. The lemma is proved. $\Box$
\begin{lemma}\label{l2} {\rm \cite[Lemma 9.13]{Fe5a})}. {\it  Let $X$ be a locally compact Abelian group containing no more than one element  of order $2$.
Let $Y$ be the character group of the group  $X$. Let $g(y)$ be a continuous function on the group $Y$  satisfying the equation
\begin{equation}\label{09_03_5n}
g(u+v)g(u-v)=g^2(u), \ u, v\in Y,
\end{equation}
and conditions
\begin{equation}\label{09_03_7n}
g(-y)=\overline{g(y)}, \ |g(y)|=1, \ y\in
Y, \ g(0)=1.
\end{equation}
 Then $g(y)$ is a character of the group $Y$.
}
\end{lemma}

We note that equation (\ref{09_03_5n}) appears in studying Gaussian distributions in the sense of Bernstein on locally compact Abelian groups (see \cite{Fe5}, \cite{Fe2012}). We also remark that generally speaking it is possible that  a continuous function $g(y)$ satisfying equation(\ref{09_03_5n}) and conditions (\ref{09_03_7n}) needs not be a character of the group $Y$, but at the same time it admits the representation
$$
g(y)=\exp\{P(y)\}, \ y\in Y,
$$
where $P(y)$ is a continuous polynomial on the group $Y$. Indeed, let $X=\mathbb{T}^2$. Then $Y\cong \mathbb{Z}^2$. Denote by $y=(m, n)$, $m, n\in \mathbb{Z}$, elements of the group $Y$. Put
$$g(y)=g(m, n)=\exp\{i\pi mn\}, \ y=(m, n)\in Y.$$

{\bf Proof of Theorem \ref{th3}}. Let $\eta_j$,  $j=1, 2, 3, 4,$ be independent random variables with values in the group   $X$ with nonvanishing characteristic functions.
 Denote by $\mu_j$ the distribution of the random variable   $\xi_j$, and by $\nu_j$ the distribution of the random variable  $\eta_j$.
 Put $M_2=b_1\eta_1+b_2\eta_2+b_3\eta_3+b_4\eta_4$,  $f_j(y)=\hat\nu_j(y)/\hat\mu_j(y)$, $\psi_j(y)=\log |f_j(y)|$, $g_j(y)=f_j(y)/|f_j(y)|$, $j=1, 2, 3, 4$.
 Obviously, each of the functions $g_j(y)$ satisfies conditions   (\ref{09_03_7n}).

 Assume that $(I)$ holds. Put $M_1=\eta_1+\eta_2+\eta_3+\eta_4$.
 Suppose that the distributions of the random vectors   $(L_1, L_2)$ and $(M_1, M_2)$ coincide.   We prove that either  $f_j(y)$ or $(f_j(y))^{-1}$ is the characteristic function of a Gaussian distribution. Thus, in the case when $(I)$ holds, the theorem will be proved.

Arguing as in the proof  of Theorem \ref{th1}, we obtain that the functions   $f_j(y)$  satisfy the equation
\begin{equation}\label{09_03_4}
f_1(u+\tilde b_1v)f_2(u+\tilde b_2v)f_3(u+\tilde b_3v)f_4(u+\tilde b_4v)=1, \ u, v\in Y.
\end{equation}
  It follows from  (\ref{09_03_4}) that the functions  $\psi_j(y)$ satisfy the equation
 \begin{equation}\label{09_03_5}
\psi_1(u+\tilde b_1v)+\psi_2(u+\tilde b_2v)+\psi_3(u+\tilde b_3v)+\psi_4(u+\tilde b_4v)=0, \ u, v\in Y.
\end{equation}
Since $Y$ is a discrete group and continuous endomorphisms  $b_j$ satisfy conditions  $(\ref{08_03_1})$, by Lemma \ref{l3}   ${(\tilde b_i-\tilde b_j)(Y)}=Y$  for $i\ne j$. By Lemma \ref{l1} $\psi_j(y)$ is a polynomial of degree at most 2. Since $\psi_j(-y)=\psi_j(y)$, $y\in
Y$, and  $\psi_j(0)=0$, it follows from the properties of polynomials, see e.g.  \cite[\S 5]{Fe5a},   that   each of the functions $\psi_j(y)$ satisfies equation    (\ref{f2}). Hence, $\psi_j(y)=\sigma_j y^2$ for some real $\sigma_j$. It means that either   $|f_j(y)|$ or  $|f_j(y)|^{-1}$  is the characteristic function of a Gaussian distribution.

We will prove that each of the functions $g_j(y)$ satisfies equation (\ref{09_03_5n}). Since any \text{\boldmath $a$}-adic solenoid $\Sigma_{{\text{\boldmath $a$}}}$  contains no more than one element of order 2, by Lemma \ref{l2} all $g_j(y)$ are characters of the group $Y$. Thus, in the case when $(I)$ holds, the theorem will be proved.

Obviously, the functions $g_j(y)$ satisfy the equation
\begin{equation}\label{09_03_6}
g_1(u+\tilde b_1v)g_2(u+\tilde b_2v)g_3(u+\tilde b_3v)g_4(u+\tilde b_4v)=1, \ u, v\in Y.
\end{equation}
Arguing as in the proof  of Theorem \ref{th1} and retaining the notation used in the proof of Theorem \ref{th1},  we eliminate from equation (\ref{09_03_6}) successively the functions $g_1(y), g_2(y), g_3(y)$.  We obtain that the function $g_4(y)$ satisfies the equation
 \begin{equation}\label{05_03_7}
{g_4(u+\tilde b_4v+(\tilde b_4-\tilde b_1)k_1+(\tilde b_4-\tilde b_2)k_2+(\tilde b_4-\tilde b_3)k_3)g_4(u+\tilde b_4v+(\tilde b_4-\tilde b_3)k_3)\over
g_4(u+\tilde b_4v+(\tilde b_4-\tilde b_2)k_2+(\tilde b_4-\tilde b_3)k_3)g_4(u+\tilde b_4v+(\tilde b_4-\tilde b_1)k_1+(\tilde b_4-\tilde b_3)k_3)}$$$$\times{g_4(u+\tilde b_4v+(\tilde b_4-\tilde b_2)k_2)g_4(u+\tilde b_4v+(\tilde b_4-\tilde b_1)k_1)\over
g_3(u+\tilde b_4v+(\tilde b_4-\tilde b_1)k_1+(\tilde b_4-\tilde b_2)k_2)g_4(u+\tilde b_4v)}=1, \ u, v\in Y.
\end{equation}
Put $v=0$, $(\tilde b_4-\tilde b_1)k_1=k$, $(\tilde b_4-\tilde b_2)k_2=l$, $(\tilde b_4-\tilde b_3)k_3=m$ in equation  (\ref{05_03_7}). Note that because
  $Y$ is a discrete group and  ${\rm Ker} (b_4-b_j)=\{0\}$,  by Lemma \ref{l3}   ${(\tilde b_4-\tilde b_j)(Y)}=Y$, $j=1, 2, 3$. Since $k_1, k_2, k_3$ are arbitrary elements
  of the group $Y$, it follows from this that the function $g_4(y)$ satisfies the equation
\begin{equation}\label{05_03_8}
{g_4(u+ k+l+m)g_4(u+m)g_4(u+l)g_4(u+k)\over
g_4(u+l+m)g_4(u+k+m)g_4(u+k+l)g_4(u)}=1, \ u, k, l, m\in Y.
\end{equation}
Put $u=-k$, $m=k$ in equation  (\ref{05_03_8}). We get
$$
{g_4(l+k)g_4(l-k)\over
g_4(l)g_4(l)}=1, \ k, l\in Y.
$$
Hence, the function $g_4(y)$ satisfies equation (\ref{09_03_5n}). For the functions $g_j(y)$, $j=1, 2, 3,$ we argue in a similar way.

Assume that $(II)$ holds. Put $M_1=\eta_1+\eta_2+\eta_3$. Arguing as in  the proof of the theorem when $(I)$ holds, we get that the characteristic function of the random vector $(L_1, L_2)$ is of the form
$$
{\bf E}[(L_1, u)(L_2, v)]=\hat\mu_1(u+\tilde b_1v)\hat\mu_2(u+\tilde b_2v)\hat\mu_3(u+\tilde b_3v)\hat\mu_3(\tilde b_4v), \ u, v\in Y,
$$
and the identical distributiveness of the random vectors    $(L_1, L_2)$ and $(M_1, M_2)$ implies that the functions $f_j(y)$ satisfy the equation
\begin{equation}\label{25_03_5}
f_1(u+\tilde b_1v)f_2(u+\tilde b_2v)f_3(u+\tilde b_3v)f_4(\tilde b_4v)=1, \ u, v\in Y.
\end{equation}
It follows from (\ref{25_03_5}) that the functions $\psi_j(y)$  satisfy the equation
 \begin{equation}\label{25_03_6}
\psi_1(u+\tilde b_1v)+\psi_2(u+\tilde b_2v)+\psi_3(u+\tilde b_3v)=-\psi_4(\tilde b_4v), \ u, v\in Y.
\end{equation}
Taking into account Lemma \ref{l3}, by Lemma \ref{l1}  (\ref{25_03_6}) implies that $\psi_j(y)$, $j=1, 2, 3$, are   polynomial   of degree at most  2. We note that since  ${\rm Ker}~b_4=\{0\}$, by Lemma \ref{l3} the set ${\tilde b_4 (Y)}$ is dense in $Y$, and hence ${\tilde b_4 (Y)}=Y$.  Taking this into account and putting in (\ref{25_03_6}) $u=0$, we obtain that $\psi_4(y)$ is also a polynomial   of degree at most  2.
As in the case when   $(I)$ holds, it follows from this that either    $|f_j(y)|$ or  $|f_j(y)|^{-1}$ is the characteristic function of a Gaussian distribution.

Obviously, the functions $g_j(y)$ satisfy the equation
\begin{equation}\label{25_03_7}
g_1(u+\tilde b_1v)g_2(u+\tilde b_2v)g_3(u+\tilde b_3v)g_4(\tilde b_4v)=1, \ u, v\in Y.
\end{equation}
We will make sure that each of the functions   $g_j(y)$ satisfies equation (\ref{09_03_5n}). Then, by Lemma \ref{l2} each of the functions $g_j(y)$ is a character of the group $Y$. Thus, in the case, when $(II)$ holds, the theorem will be also proved. First verify that the function  $g_4(y)$ satisfies equation (\ref{09_03_5n}). Arguing as in the proof  of the theorem when $(I)$ holds,  we eliminate from equation  (\ref{25_03_7}) the functions $g_1(y)$ and  $g_2(y)$. We get
\begin{equation}\label{25_03_8}
{g_3(u+\tilde b_3v+(\tilde b_3-\tilde b_1)k_1+(\tilde b_3-\tilde b_2)k_2)g_4(\tilde b_4v+\tilde b_4k_1+\tilde b_4k_2)g_3(u+\tilde b_3v)g_4(\tilde b_4v)\over
g_3(u+\tilde b_3v+(\tilde b_3-\tilde b_2)k_2)g_4(\tilde b_4v+\tilde b_4k_2)g_3(u+\tilde b_3v+(\tilde b_3-\tilde b_1)k_1)g_4(\tilde b_4v+\tilde b_4k_1)}=1, \ \ u, v\in Y.
\end{equation}
Eliminating from equation  (\ref{25_03_8}) the function $g_3(y)$, we obtain
\begin{equation}\label{25_03_9}
{g_4(\tilde b_4v+\tilde b_4k_1+\tilde b_4k_2+\tilde b_4k_3)g_4(\tilde b_4v+\tilde b_4k_3)g_4(\tilde b_4v+\tilde b_4k_2)g_4(\tilde b_4v+\tilde b_4k_1)\over
g_4(\tilde b_4v+\tilde b_4k_2+\tilde b_4k_3)g_4(\tilde b_4v+\tilde b_4k_1+\tilde b_4k_3)g_4(\tilde b_4v+\tilde b_4k_1+\tilde b_4k_2)g_4(\tilde b_4v)}=1, \ u, v\in Y.
\end{equation}
Note that since ${\rm Ker}~b_4=\{0\}$, then by Lemma \ref{l3} the set ${\tilde b_4 (Y)}$ is dense in $Y$, and hence ${\tilde b_4 (Y)}=Y$. Put in   (\ref{25_03_9}) $\tilde b_4v=u$, $\tilde b_4k_1=k$, $\tilde b_4k_2=l$, $\tilde b_4k_3=m$. Since $v, k_1, k_2, k_3$ are arbitrary elements
  of the group $Y$, we  get that the function  $g_4(y)$ satisfies equation (\ref{05_03_8}). The desired assertion follows from this.

Make sure now that the functions $g_j(y)$, $j=1, 2, 3$, also satisfy equation (\ref{09_03_5n}). Let $k$ be an arbitrary element of the group
$Y$.   Substitute $u+k$ for $u$  in equation (\ref{25_03_8}).
Dividing the resulting equation into equation (\ref{25_03_8}), we find
\begin{equation}\label{25_03_10}
{g_3(u+k+\tilde b_3v+(\tilde b_3-\tilde b_1)k_1+(\tilde b_3-\tilde b_2)k_2)g_3(u+k+\tilde b_3v)\over
g_3(u+k+\tilde b_3v+(\tilde b_3-\tilde b_2)k_2)g_3(u+k+\tilde b_3v+(\tilde b_3-\tilde b_1)k_1)}$$$$\times{g_3(u+\tilde b_3v+(\tilde b_3-\tilde b_2)k_2)g_3(u+\tilde b_3v+(\tilde b_3-\tilde b_1)k_1)\over g_3(u+\tilde b_3v+(\tilde b_3-\tilde b_1)k_1+(\tilde b_3-\tilde b_2)k_2)g_3(u+\tilde b_3v)}=1, \ \ u, v\in Y.
\end{equation}
Put in (\ref{25_03_10}) $v=0$, $(\tilde b_3-\tilde b_1)k_1=l$, $(\tilde b_3-\tilde b_2)k_2=m$. Taking into account Lemma \ref{l3} and conditions (\ref{25_03_2}), the resulting equation implies that the function  $g_3(y)$ satisfies equation (\ref{05_03_8}).  The desired assertion follows from this. For the functions $g_j(y)$, $j=1, 2$, we argue in a similar way.
The theorem is completely proved. $\Box$

\begin{remark}\label{re5}.
{\rm The reasoning  similar to that in   Remark \ref{re3} shows that
{\rm Theorem \ref{th3}    fails  if we omit conditions    (\ref{25_03_1}) and  (\ref{25_03_2}).}}
\end{remark}

\begin{remark}\label{re6}.
{\rm Generally speaking,
{\rm Theorem \ref{th3}    fails  for an arbitrary locally compact Abelian group. Let $X=\mathbb{R}^2$. Then $Y\cong X$. Denote by $x=(x_1, x_2)$, $y=(y_1, y_2)$ elements of the groups $X$ and $Y$ respectively. Put $b_j = \left(\begin{matrix}j&0\\ 0&-j\end{matrix}\right)$, j=1, 2, 3, 4. Then $b_j$ are continuous endomorphisms of the group $X$, satisfying conditions (\ref{25_03_1}). Put $\mu_1=\mu_2=\mu_3=\mu_4=\mu$, where $\mu$ is a distribution on the group $X$ with the characteristic function $\hat\mu(y_1, y_2)=\exp\{-4(y_1^2+y_2^2)\}$. Let $\nu_j$ be distributions on the group $X$ with the characteristic functions $\hat\nu_1(y_1, y_2)=\exp\{-3y_1^2-5y_2^2\}$, $\hat\nu_2(y_1, y_2)=\exp\{-7y_1^2-y_2^2\}$, $\hat\nu_3(y_1, y_2)=\exp\{-y_1^2-7y_2^2\}$, $\hat\nu_4(y_1, y_2)=\exp\{-5y_1^2-3y_2^2\}$. Let $\xi_j$ be independent random variables with values in the group   $X$ and distributions $\mu_j$, and $\eta_j$ be independent random variables with values in the group   $X$ and distributions $\nu_j$, $j=1, 2, 3, 4.$ Put $L_1=\xi_1+\xi_2+\xi_3+\xi_4$, $L_2=b_1\xi_1+b_2\xi_2+b_3\xi_3+b_4\xi_4$, $M_1=\eta_1+\eta_2+\eta_3+\eta_4$, $M_2=b_1\eta_1+b_2\eta_2+b_3\eta_3+b_4\eta_4$. On the one hand, it is easy to verify that (\ref{05_03_2}) holds. Hence the random vectors  $(L_1, L_2)$ and $(M_1, M_2)$ are identically distributed. On the other hand, it is obvious that do not exist Gaussian distributions $\gamma_j$ on $X$ such that either $\mu_j=\nu_j*\gamma_j$ or $\nu_j=\mu_j*\gamma_j$.}}
\end{remark}

\bigskip

\centerline{\bf 4. Comments}

\bigskip

In the article \cite{PR} B.L.S.~Prakasa Rao noted that the proof by I.~Kotlarski of Theorem B can be extended without any changes to locally compact Abelian groups. Our proof of Theorem \ref{th1} is new for the group of real numbers $\mathbb{R}$ and differs from the proof of Theorem A given by C.R.~Rao in \cite{CRR}. In \cite{CRR} for solving equations  (\ref{05_03_3}) and (\ref{06_03_1})  one takes the logarithm of both sides of the equation considering the main branch of it. Then one proves that $g_k(y)=\log f_k(y)$ are linear functions, and hence $f_k(y)=e^{ix_ky}$, $x_k\in \mathbb{R}$, $k=1, 2, 3$. Theorem A follows from this. In the case of locally compact Abelian groups, generally speaking, one can not take the logarithm of equations  (\ref{05_03_3}) and (\ref{06_03_1}). Moreover, since $X$ is an arbitrary locally compact Abelian group, it is well-known that, generally speaking,  not every character $(x, y)$ of the group $Y$ is of the form \begin{equation}\label{16_07_2}
(x, y)=e^{ih(y)}, \ y\in Y,
\end{equation}
where $h(y)$ is a continuous real valued function satisfying the equation \begin{equation}\label{16_07_1}
h(u+v)=h(u)+h(v), \ u, v\in Y.
\end{equation}

We emphasize the following. In    \cite{CRR} C.R.~Rao  proves that if $\xi_j$, $j=1, 2, \dots, n$, are independent random variables with nonvanishing characteristic functions, then the distributions of two linear forms of $\xi_j$ determine  the characteristic functions of these random variables  up to a factor of the form $\exp\{P_j(s)\}$, where $P_j(s)$ is a polynomial of degree   at most $n-2$. In particular, if we have three independent random variables, then the degree of $P_j(s)$ is at most 1. This easily implies that $P_j(s)=ib_js$, $b_j\in \mathbb{R}$. Hence the distributions  of the random variables  are determined  up to a  change of location. If we have four independent random variables,  then the degree of $P_j(s)$ is at most 2, and this   implies that  $P_j(s)=a_js^2+ib_js$, $a_j, b_j\in \mathbb{R}$. Hence the distributions  of the random variables  are determined  up to a convolution with a Gaussian distribution.  It is interesting to note that, in spite of the fact that   Theorem \ref{th1} is a complete analogue of the corresponding statement on the real line, generally speaking, for an arbitrary locally compact Abelian group  $X$ it is not true that the distributions of two linear forms of three independent random variables with nonvanishing characteristic functions determine  the characteristic functions of these random variables up to a factor of the form   $\exp\{P(y)\}$, where $P(y)$ is a polynomial  on the group $Y$. It is easy to construct the corresponding example on the    \text{\boldmath $a$}-adic solenoid $X=\Sigma_{{\text{\boldmath $a$}}}$ for which its character group $Y=H_{\text{\boldmath $a$}}$ is topologically isomorphic to the discrete group of all rational numbers. This is due to the fact that for an arbitrary    \text{\boldmath $a$}-adic solenoid $X=\Sigma_{{\text{\boldmath $a$}}}$, there is its character $(x, y)$ which is not represented in the form (\ref{16_07_2}), where   $h(y)$ is a continuous real valued function satisfying the equation (\ref{16_07_1}).

\bigskip

I want to thank the reviewer for a very careful reading of the article and useful remarks.

\end{document}